\documentclass[12pt]{amsart}
\usepackage{amssymb}
\usepackage{amsmath}
\usepackage{mathabx}
\usepackage{amsthm}
\usepackage{amstext}
\usepackage{amsopn}
\usepackage[dvips]{graphicx}
\usepackage{mathrsfs}
\usepackage{latexsym}
\DeclareMathAlphabet{\mathpzc}{OT1}{pzc}{m}{it}
\allowdisplaybreaks
\newtheorem{Definition}{Definition}[section]
\newtheorem{Proposition}{Proposition}[section]
\newtheorem{Lemma}{Lemma}[section]
\newtheorem{Theorem}{Theorem}[section]
\newtheorem{Corollary}{Corollary}[section]
\newtheorem{Remark}{Remark}[section]
\newtheorem{Example}{Example}[section]

\begin{document}
\bibliographystyle{plain}
\title[Multivariate Fuss-Narayana polynomials]
{Multivariate Fuss-Narayana polynomials and \\their application to random matrices}
\author[R. Lenczewski, R. Sa{\l}apata]{Romuald Lenczewski and Rafa{\l} Sa{\l}apata}
\address{Romuald Lenczewski \newline
Instytut Matematyki i Informatyki, Politechnika Wroc\l{}awska, \newline
Wybrze\.{z}e Wyspia\'{n}skiego 27, 50-370 Wroc{\l}aw, Poland 
\newline
\indent
Rafa{\l} Sa{\l}apata \newline
Instytut Matematyki i Informatyki, Politechnika Wroc\l{}awska, \newline
Wybrze\.{z}e Wyspia\'{n}skiego 27, 50-370 Wroc{\l}aw, Poland  \vspace{10pt}}
\email{Romuald.Lenczewski@pwr.wroc.pl, Rafal.Salapata@pwr.wroc.pl}
\begin{abstract}
It has been shown recently that the limit moments of $W(n)=B(n)B^{*}(n)$, where $B(n)$ is 
a product of $p$ independent rectangular random matrices, are certain homogenous polynomials 
$P_{k}(d_0,d_1, \ldots , d_{p})$ in the asymptotic dimensions of these matrices. 
Using the combinatorics of noncrossing partitions, we explicitly determine these polynomials and show that they are
closely related to polynomials which can be viewed as {\it multivariate Fuss-Narayana polynomials}.
Using this result, we compute the moments of $\varrho_{t_1}\boxtimes \varrho_{t_2}\boxtimes \ldots \boxtimes \varrho_{t_p}$
for any positive $t_1,t_2, \ldots , t_n$, where $\boxtimes$ is the free multiplicative convolution in free probability and $\varrho_{t}$ is
the Marchenko-Pastur distribution with shape parameter $t$.
\end{abstract}

\maketitle
\section{Introduction}
The motivation of this paper comes from some recent developments concerning random matrices and their products [8], where
combinatorial formulas for the polynomials studied in this paper were obtained.

For any given $p\in {\mathbb N}$ and any $n\in {\mathbb N}$, 
consider the product of independent rectangular Gaussian random matrices
$$
B(n)=X_1(n)X_2(n)\ldots X_{p}(n)
$$
where $n\in {\mathbb N}$ and their dimensions are such that the product is well-defined. 
If $X_{j}$ is an $N_{j-1}(n) \times N_{j}(n)$ matrix for any $j$, we assume that 
$\lim_{n\rightarrow \infty} N_{j}(n)/n =d_j> 0$, where $d_0, d_1, \ldots d_{p}$ are called {\it asymptotic dimensions}. 
This notation slightly differs from $d_1,d_2, \ldots , d_{p+1}$ used in [8], but it is more convenient 
in combinatorial formulas.
Finally, let $\tau_{0}(n)$ be the trace over the set of $N_{0}(n)$ basis vectors composed with classical expectation. 

It is shown in [8] that under certain natural assumptions 
$$
\lim_{n \rightarrow \infty}\tau_{0}(n)\left(\left(B(n)B^{*}(n)\right)^{k}\right)=P_{k}(d_0,d_1, \ldots , d_{p})
$$
where on the right-hand side we have a certain homogenous polynomial of order $2pk$ for any
$k\in \mathbb{N}$. There, instead of the matrices $X_1(n), X_2(n),\ldots ,X_{p}(n)$, 
symmetric blocks embedded in a large $n\times n$ square matrix were used and $\tau_{0}(n)$ was the  
{\it partial trace} over a subset of basis vectors composed with classcial expectation. 
That formulation followed from the fact 
that it was the ensemble of symmetric blocks which was shown to converge in moments 
to the ensemble of operators generalizing free Gaussian operators
called {\it matricially free Gaussian operators}.

Non-homogenous polynomials obtained from $P_{k}$ by dividing them by $d_{0}^{kp}$
were called `generalizations of Narayana polynomials' since for
$p=1$ they become the well-known Narayana polynomials of one variable corresponding 
to Catalan numbers and their decompositions in terms of Narayana numbers. 
In random matrix theory, these well-known combinatorial 
objects are important since they are related to the moments of the Marchenko-Pastur distribution, 
the limit distribution of Wishart random matrices [10].

In the case when $p>1$ and all matrices in the product are square, the limit moments
are Fuss-Narayana polynomials of one variable. These polynomials correspond, in turn,
to Fuss-Catalan numbers and their decompositions in terms of Fuss-Narayana numbers. 
For arbitrary asymptotic dimensions, it is therefore natural to expect that we should obtain some
multivariate analogs of Fuss-Narayana polynomials.

A combinatorial definition of the polynomials $P_{k}$ was given in [8]. In this paper, we determine their explicit 
form, using purely combinatorial methods. Namely, we show that
$$
P_{k}(d_0,d_{1},\ldots , d_{p})=\sum_{j_0 + \ldots + j_{p} = pk+1}\frac{1}{k}{k \choose j_0}{k \choose j_1}\ldots {k \choose j_{p}}
\;d_{0}^{j_{0}-1}d_{1}^{j_1}\ldots d_{p}^{j_{p}}
$$
where the indices $j_0, j_1, \ldots , j_p$ are natural numbers. When we divide these polynomials by $d_{0}^{kp}$, which 
corresponds to a different normalization of random matrices, we obtain 
$$
F_{k}(t_{1},t_{2}, \ldots , t_{p})=d_{0}^{-kp}P_{k}(d_{0},d_{1}, \ldots , d_{p}),
$$
where $t_j=d_{j}/d_{0}$ for any $j$, called {\it multivariate Fuss-Narayana polynomials} 
since for $p=1$ they become the well-known Narayana polynomials.

The coefficients built from binomial expressions play the role of {\it generalized Fuss-Narayana numbers}.
Such numbers appeared recently in a different context in the paper of Loktev on Weyl modules of Lie algebras [9], who showed that they were dimensions of weight spaces of two-variable Weyl modules of Lie algebras $gl_{p+1}$.

In random matrix theory, Fuss-Catalan numbers appeared
in the paper of Alexeev, Goetze and Tikhomirov [1], who showed, using the methods of classical probability, 
that they were the limit moments of the asymptotic distribution of squared singular values
of powers of a random matrix with independent entries. Moreover, similar techniques 
were used by the same authors to study the case of products of independent random matrices [2]. 
In the general case of rectangular matrices, the limit distribution for singular values of such products 
has only been described by an algebraic equation for its Cauchy-Stieltjes transform and a formula for 
moments has not been derived.  
Our approach is quite different since we use noncommutative probability and operator algebras. 
This allows us to formulate the problem in a purely combinatorial fashion and then 
solve it using combinatorial techniques. 

An important probabilistic context in which Fuss-Narayana polynomials of one variable appear is that of free Bessel laws 
$$
\pi_{p,\,t}=\varrho_{1}^{\,\boxtimes (p-1)}\boxtimes \varrho_{t},
$$
where $p\in {\mathbb N}$, $\varrho_{t}$ is the Marchenko-Pastur distribution with shape parameter $t>0$,
and $\boxtimes$ denotes the free multiplicative convolution.
They were defined by Banica, Belinschi, Capitaine and Collins [3], who have shown that their moments are given by Fuss-Narayana polynomials in $t$. 

A more general multiplicative free convolution of Marchenko-Pastur laws is of the form
$$
\varrho_{t_1}\boxtimes \varrho_{t_2}\boxtimes \ldots \boxtimes \varrho_{t_p}
$$
for any positive $t_1,t_2,\ldots , t_p$. We show in this paper that the moments of such convolutions
are given by multivariate Fuss-Narayana polynomials. For the foundations of free probability and, in particular,
for the definition of the free multiplicative convolution, see [11].

\section{Multivariate Fuss-Narayana polynomials}
The {\it Fuss-Catalan numbers} associated with $p\in {\mathbb N}$ are given by the formula 
$$
C_k = \frac{1}{k} {(p+1)k \choose pk+1}
$$
where $k\in {\mathbb N}$. The following decomposition of Catalan numbers is known as the generalized Vandermonde's identity.

\begin{Proposition}
For any $p,k\in {\mathbb N}$, the following decomposition hold:
$$
C_{k}= \sum_{j_0 + \ldots + j_p = pk+1}
N(k,j_0, \ldots , j_p)
$$
where the summation runs over all $j_0,\ldots , j_p\in [k]=\{1,2, \ldots , k\}$, for which it holds that
$j_0+\ldots +j_p=pk+1$ and for such values
$$
N(k,j_0, \ldots , j_p)=\frac{1}{k} {k \choose j_0}{k \choose j_1}\ldots {k\choose j_p}.
$$
These numbers will be called generalized Fuss-Narayana numbers.
\end{Proposition}
{\it Proof.}
If we apply the Vandermonde's identity to the formula for the Fuss-Catalan numbers $p$ times, we obtain the above formula.
\hfill $\blacksquare$\\

Note that after the first application of the Vandermonde's identity, we get the decomposition
$$
C_k= 
\sum_{j_0 +j_1 = pk+1}
\frac{1}{k} 
{k \choose j_{0}}{pk \choose j_1},
$$
where the summands on the right-hand side are called {\it Fuss-Narayana numbers}.

It will be convenient to use vector notations 
$$
{\bf d} = (d_0,d_1,\ldots,d_p),\;{\bf j} = (j_0,j_1,\ldots,j_p)
$$
where $d_0,d_1, \ldots , d_p$ are variables and $j_0,j_1, \ldots ,j_p$ nonnegative integers. 
Using these vectors to write the generalized Fuss-Narayana numbers as 
$N(k,{\bf j})$ and generalized powers in the form
$$
{\bf d^j} = d_0^{j_0} d_1^{j_1} \ldots d_p^{j_p},
$$
we can define the generating function for the generalized Fuss-Narayana numbers
by 
$$
\mathcal{N}(x,{\bf d}) = \sum_{k=1}^{\infty}
\sum_{j_0 + \ldots + j_p = pk+1}
N(k,{\bf j})\,{\bf d^j}x^{k}.
$$
We will set $N(k,{\bf j})=0$ whenever $j_0 + \ldots + j_p \neq pk+1$ or $j_i\notin [k]$ for some $i$.

Dividing the homogenous polynomials in $p+1$ asymptotic dimensions of the form
$$
R_{k}(d_0, d_1,\ldots , d_{p})=\sum_{j_0 + \ldots + j_p = pk+1}
N(k,{\bf j})\,d_{0}^{j_{0}}d_{1}^{j_{1}}\ldots d_{p}^{j_{p}}
$$
by $d_{0}^{pk+1}$, we obtain certain non-homogenous polynomials in $p$ variables $t_{i}=d_{i}/d_{0}$, where $i\in [p]$, 
which play the role of mutlivariate generalizations of Narayana polynomials, obtained for $p=1$. 
\begin{Definition}
{\rm Non-homogenous polynomials of $p$ variables of the form
$$
F_{k}(t_1,t_2, \ldots , t_p)= \sum_{j_0 + \ldots + j_p = pk+1}N(k,j_{0}, \ldots , j_{p})t_{1}^{j_1}t_{2}^{j_2}\ldots t_{p}^{j_p},
$$
where $k\in {\mathbb N}$, will be called {\it multivariate Fuss-Narayana polynomials}. }
\end{Definition}

\begin{Remark}
{\rm Of course, one can also obtain these polynomials from $R_k$ by setting $d_{0}=1$ and $d_j=t_p$ for $j\in [p]$. In fact, 
there is a number of equivalent definitions of this type, namely
$$
F_{k}(t_1,t_2, \ldots , t_p)=R_{k}(1,t_{1}, \ldots , t_{p})=\ldots =R_{k}(t_{1},t_{2}, \ldots , t_{p},1)
$$
due to the symmetric form of generalized Fuss-Narayana numbers.
}
\end{Remark}
In order to establish a formula satisfied by the generating function $\mathcal{N}$, we will use the next theorem, which is a 
special version of the well known Lagrange Inversion Theorem 
(see, for instance, Appendix A in the paper of Deutch [4] or the book of Wilf [12]).

\begin{Theorem} 
Assume that a generating function $f(z)$ satisfies the functional equation 
\begin{equation}
f(z) = zH(f(z))\,,
\end{equation}
where $H$ is a polynomial. Then $(1)$ has a unique solution and
$$
[z^n] f(z) = \frac{1}{n} [\lambda^{n-1}] H^n(\lambda)\,,
$$
where $[z^n] f(z)$ is the coefficient of the series $f(z)$ standing by $z^n$.
\end{Theorem}

\begin{Proposition}
For any $d_0, \ldots , d_p$, the unique solution of the equation
\begin{equation}
g(x) = x\prod_{i=0}^{p}(g(x) + d_i)
\end{equation}
is given by $g(x)=\mathcal{N}(x,{\bf d})$.
\end{Proposition}
{\it Proof.} Applying Theorem 2.1 to the polynomial $H(\lambda) = (\lambda + d_0)(\lambda + d_1)\ldots(\lambda + d_p)$, we get existence of a unique solution $g(x)$ of equation (2). 
We must have
\begin{eqnarray*}
[x^n] g(x) &=& \frac{1}{n} [\lambda^{n-1}] 
\left(\prod_{i=0}^{p}\left(\lambda + d_i\right)\right)^{n}\\
&=& \frac{1}{n} [\lambda^{n-1}] \prod_{i=0}^{p}\left(\sum_{k=0}^{n} {n \choose k} \lambda^{k} d_i^{\,n-k} \right) \\
&=& \frac{1}{n} [\lambda^{n-1}]  \sum_{k=0}^{\infty} \left( \sum_{k_0 + \ldots + k_p = k} 
{n \choose k_0}\ldots {n\choose k_p} d_0^{\,n-k_0} \ldots d_p^{\,n-k_p}  \right) \lambda^{k}  \\
&=& \frac{1}{n}  \sum_{j_0 + \ldots + j_p = pn+1} 
{n \choose j_0} \ldots {n \choose j_p} d_0^{j_0} \ldots d_p^{j_p},
\end{eqnarray*}
where we put $n-k_i = j_i$ for $i\in\{0,1,\ldots,p\}$ to get the last equation. This proves
that $g(x)={\mathcal N}(x, {\bf d})$. \hfill $\blacksquare$\\

\section{Lemmas}
In this section we shall prove two combinatorial lemmas, in which we compare the generating function ${\mathcal N}$ of Section 2 defined by polynomials $R_k$
with the generating function ${\mathcal N}_{0}$ defined by polynomials $P_k$.

By a noncrossing pair partition of the set $[m]:=\{1,2, \ldots , m\}$, where $m=2k$ is an even natural number,  
we understand a collection $\pi=\{\pi_1, \pi_2, \ldots , \pi_k\}$ of disjoint two-element subsets of $[m]$ called blocks, 
such that there are no blocks $\{i,j\}$ and $\{p,q\}$ for which $i<p<j<q$.
The set of such partitions will be denoted by $\mathcal{NC}_{m}^{2}$.
If $\{i,j\}$ is a block and $i<j$, then $i$ and $j$ are called the left and the right legs of this block, respectively.
For a given noncrossing pair partition $\pi$, we denote by $\mathcal{R}(\pi)$ the set of its right legs.

We will consider certain noncrossing pair-partitions of $[m]$, where $m=2pk$ and $p, k\in {\mathbb N}$, which are associated with 
words of the form
$$
W_0^{k}=(1\ldots pp^{*}\ldots 1^{*})^{k}
$$
built from $p$ starred and $p$ unstarred letters. More generally, we will denote by $W_i$ the word that arises from $W_0$ by the cyclic shift of its letters to 
the right by $i$ positions, where $i\in \{0,1,\ldots, p\}$, namely
$$
W_{i}=i^{*}\ldots 1^{*}1 \ldots pp^{*}\ldots (i+1)^{*},
$$
thus, in particular, $W_p =p^*  \ldots 1^* 1 \ldots p$. We will also consider powers $W_{i}^{k}$ of such shifted words. 
For simplicity, we supress $p$ in all these notations, in contrast to the notation used in [8].

\begin{Definition}
{\rm We shall say that $\pi \in \mathcal{NC}^2_{2pk}$ is {\it adapted} to the word $W_i^k$ if all its blocks 
are associated with pairs of letters of the form $\{l,l^*\}$ for some $l\in \{1,\ldots,p\}$. 
By $\mathcal{NC}^2_{2pk}(W^k_i)$ we denote the set of all pair partitions from $\mathcal{NC}^2_{2pk}$ which are
adapted to $W^k_i$. We set $\mathcal{NC}^2_{0}(W^0_i)=\{\emptyset\}$.}
\end{Definition}

Thus, if $\pi\in \mathcal{NC}^2_{2pk}(W^k_i)$, then blocks of $\pi$ are pairs $\{r,s\}$ 
in which $r$ is associated with the letter $l$ if and only if $s$ is associated with the letter $l^{*}$. 
Thus, if $r<s$, then $s$ is the right leg of this block and it is associated with $l^{*}$, whereas if 
$r>s$, then $r$ is the right leg of this block and it is associated with $l$. Therefore, it is meanigful to 
define the sets of right legs of $\pi$ associated with $l$ and $l^{*}$, respectively, and denote these sets
by ${\mathcal R}_{l}(\pi)$ and ${\mathcal R}_{l}^{*}(\pi)$. Using these sets, whose union gives ${\mathcal R}(\pi)$, we 
defined in [8] the following family of homogenous polynomials in $p+1$ variables.

\begin{Definition}
{\rm Define polynomials in variables $d_0,d_1, \ldots , d_{p}$ of the form
$$
P_{k}(d_0,d_1, \ldots , d_{p})=\sum_{\pi \in \mathcal{NC}_{2kp}^{2}(W_{0}^{k})}d_{0}^{\,j_0(\pi)}d_{1}^{\,j_1(\pi)}\ldots d_{p}^{\,j_{p}(\pi)}
$$
for any $k,p\in \mathbb{N}$, where
$$
j_l(\pi)=|\mathcal{R}_{l+1}(\pi)|+|\mathcal{R}_{l}^{*}(\pi)|
$$
for any $\pi\in \mathcal{NC}_{2kp}^{2}(W_{0}^{k})$ and $0\leq l\leq p$, where we set $\mathcal{R}_{0}^{*}(\pi)=\emptyset$ and 
$\mathcal{R}_{p+1}(\pi)=\emptyset$. }
\end{Definition}

Our goal is now to count the number of noncrossing pair partitions $\pi \in \mathcal{NC}_{2kp}^{2}(W_{0}^{k})$ 
which contribute identical monomials to the polynomials $P_k$. In other words, if
we label the legs of $\pi$ by letters from the set $\{1,\ldots , p, p^{*}, \ldots ,1^{*}\}$, we would like to count 
the blocks labelled by ordered pairs $(l,l^*)$ or $((l+1)^*,(l+1))$ for any $l\in\{0,1,\ldots,p\}$. 

We are mainly interested in counting the noncrossing pair partitions adapted to $W_{0}^{k}$ which contribute identical monomials, but
we shall need a family which encodes detailed information about the labellings of their right legs, namely
$$
N_i(k,{\bf j}) = \#\{ \pi \in \mathcal{NC}^2_{2pk}(W_i^k); j_l(\pi) = j_l \ \text{for} \ l\in \{0,1,\ldots,p\}\}
$$
where ${\bf j}=(j_0, j_1, \ldots , j_p)$ and $i \in \{0,1,\ldots,p\}$. It is obvious that 
$$
N_i(k,{\bf j})=0\;\;{\rm whenever}\;\;j_0 + \ldots + j_p \neq pk\;\;{\rm or}\;\;j_{i}> k\;\;{\rm for}\;{\rm some}\;i,
$$ 
and that $N_i(0,0,\ldots,0)=1$ for any $i$. Let us add that the reason why we count together 
right legs labelled by $l^{*}$ and $l+1$ follows from the way we multiply 
rectangular matrices and their adjoints and is clear from the proof of [8, Theorem 10.1].

The corresponding generating functions are given by
$$
\mathcal{N}_i(x,{\bf d}) = \sum_{k,j_0,\ldots,j_p=0}^{\infty} N_i(k,{\bf j}) {\bf d^j}x^{k}
$$
We are especially interested in the coefficients of the generating function for $i=0$, 
i.e. we will study homogenous polynomials of degrees $pk$ in variables $d_0, d_1,\ldots, d_p$ given by
the combinatorial formulas 
\begin{equation}
P_k(d_0,d_1,\ldots , d_p) = \sum_{\pi\in\mathcal{NC}^2_{2pk}(W_0^k)}   d_0^{j_0(\pi)} d_1^{j_1(\pi)} \ldots d_{p}^{j_{p}(\pi)}  = 
\sum_{j_0,\ldots,j_p=0}^{\infty} N_0(k,{\bf j}) {\bf d^j}
\end{equation}
for $k\geqslant0$, where, in the last formula, only finitely many terms do not vanish. Obviously, 
$$
\mathcal{N}_0(x,{\bf d}) = \sum_{k=0}^{\infty} P_k(d_0,d_1, \ldots , d_p) x^k.
$$
Our goal in this section is to find the coefficients $N_0(k,{\bf j})$ of the polynomials $P_k$ and show 
that they are generalized Fuss-Narayana numbers. Moreover, if we divide $P_k$ by $d_{0}^{kp}$, 
we will obtain multivariate Fuss-Narayana polynomials of $t_1=d_1/d_0, \ldots , t_p=d_p/d_0$.

\begin{Example}
{\rm Consider the case $p=2$. Then, there are three noncrossing pair partitions which are adapted to the word $W_{0}^{2}=(122^{*}1^{*})^{2}$, 
corresponding to the diagrams

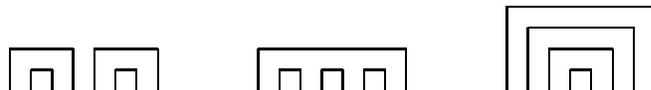
\begin{picture}(250.00,55.00)(00.00,5.00)
\put(76.00,10.00){\line(0,1){16.00}}
\put(84.00,10.00){\line(0,1){8.00}}
\put(92.00,10.00){\line(0,1){8.00}}
\put(100.00,10.00){\line(0,1){16.00}}

\put(76.00,26.00){\line(1,0){24.00}}
\put(84.00,18.00){\line(1,0){8.00}}

\put(108.00,10.00){\line(0,1){16.00}}
\put(116.00,10.00){\line(0,1){8.00}}
\put(124.00,10.00){\line(0,1){8.00}}
\put(132.00,10.00){\line(0,1){16.00}}

\put(108.00,26.00){\line(1,0){24.00}}
\put(116.00,18.00){\line(1,0){8.00}}

\put(170.00,10.00){\line(0,1){16.00}}
\put(178.00,10.00){\line(0,1){8.00}}
\put(186.00,10.00){\line(0,1){8.00}}
\put(194.00,10.00){\line(0,1){8.00}}
\put(202.00,10.00){\line(0,1){8.00}}
\put(210.00,10.00){\line(0,1){8.00}}
\put(218.00,10.00){\line(0,1){8.00}}
\put(226.00,10.00){\line(0,1){16.00}}

\put(170.00,26.00){\line(1,0){56.00}}
\put(178.00,18.00){\line(1,0){8.00}}
\put(194.00,18.00){\line(1,0){8.00}}
\put(210.00,18.00){\line(1,0){8.00}}

\put(264.00,10.00){\line(0,1){32.00}}
\put(272.00,10.00){\line(0,1){24.00}}
\put(280.00,10.00){\line(0,1){16.00}}
\put(288.00,10.00){\line(0,1){8.00}}
\put(296.00,10.00){\line(0,1){8.00}}
\put(304.00,10.00){\line(0,1){16.00}}
\put(312.00,10.00){\line(0,1){24.00}}
\put(320.00,10.00){\line(0,1){32.00}}

\put(264.00,42.00){\line(1,0){56.00}}
\put(272.00,34.00){\line(1,0){40.00}}
\put(280.00,26.00){\line(1,0){24.00}}
\put(288.00,18.00){\line(1,0){8.00}}

\end{picture}
\\
and thus
$$
P_{2}(d_0,d_1,d_2)=d_{1}^{\,2}d_{2}^{\,2}+d_{0}d_{1}d_{2}^{\,2}+ d_0d_{1}^{\,2}d_{2}
$$
is the associated polynomial. If we take $W_{0}^{3}=(122^{*}1^{*})^{3}$, then 
the number of adapted noncrossing pair partitions grows to 12, which will coincide with 
the corresponding Fuss-Catalan number of Corollary 4.2. Let us just give the corresponding polynomial
$$
P_3(d_0,d_1,d_2) = d_1^3 d_2^3 + 3 d_0 d_1^2 d_2^3 + 3 d_0 d_1^3 d_2^2 + d_0^2 d_1 d_2^3 + 3 d_0^2 d_1^2 d_2^2 + d_0^2 d_1^3 d_2
$$
since drawing all the diagrams seems too elaborate.
}
\end{Example}

For further purposes, we introduce a function $\varphi: \mathcal{NC}^2_{2pk} \rightarrow \mathcal{NC}^2_{2pk}$ defined by the diagram

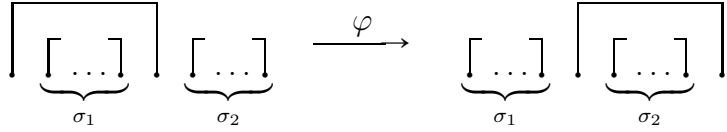
\unitlength=0.8mm
\special{em:linewidth 0.4pt}
\linethickness{0.4pt}
\begin{picture}(120.00,37.00)(35.00,26.00)

\put(127.00,45.00){$\rightarrow$}
\put(116.00,46.20){\line(1,0){12.00}}
\put(122.50,48.50){$\varphi$}

\put(66.00,41.00){\line(0,1){12.00}}
\put(66.00,53.00){\line(1,0){24.00}}
\put(90.00,41.00){\line(0,1){12.00}}
\put(66.00,41.00){\circle*{1.00}}
\put(90.00,41.00){\circle*{1.00}}

\put(72.00,41.00){\line(0,1){6.00}}
\put(72.00,47.00){\line(1,0){2.00}}
\put(84.00,41.00){\line(0,1){6.00}}
\put(84.00,47.00){\line(-1,0){2.00}}
\put(84.00,41.00){\circle*{1.00}}
\put(72.00,41.00){\circle*{1.00}}

\put(96.00,41.00){\line(0,1){6.00}}
\put(96.00,47.00){\line(1,0){2.00}}
\put(96.00,41.00){\circle*{1.00}}
\put(108.00,41.00){\line(0,1){6.00}}
\put(108.00,47.00){\line(-1,0){2.00}}
\put(108.00,41.00){\circle*{1.00}}

\put(75.70,41.00){$\ldots$}
\put(99.50,41.00){$\ldots$}
\put(70.30,40.00){$\underbrace{\hspace{35pt}}_{\sigma_1}$}
\put(94.30,40.00){$\underbrace{\hspace{35pt}}_{\sigma_2}$}


\put(184.00,41.00){\line(0,1){12.00}}
\put(184.00,53.00){\line(-1,0){24.00}}
\put(160.00,41.00){\line(0,1){12.00}}
\put(184.00,41.00){\circle*{1.00}}
\put(160.00,41.00){\circle*{1.00}}

\put(142.00,41.00){\line(0,1){6.00}}
\put(142.00,47.00){\line(1,0){2.00}}
\put(154.00,41.00){\line(0,1){6.00}}
\put(154.00,47.00){\line(-1,0){2.00}}
\put(154.00,41.00){\circle*{1.00}}
\put(142.00,41.00){\circle*{1.00}}

\put(166.00,41.00){\line(0,1){6.00}}
\put(166.00,47.00){\line(1,0){2.00}}
\put(166.00,41.00){\circle*{1.00}}
\put(178.00,41.00){\line(0,1){6.00}}
\put(178.00,47.00){\line(-1,0){2.00}}
\put(178.00,41.00){\circle*{1.00}}

\put(145.70,41.00){$\ldots$}
\put(169.50,41.00){$\ldots$}
\put(140.30,40.00){$\underbrace{\hspace{35pt}}_{\sigma_1}$}
\put(164.30,40.00){$\underbrace{\hspace{35pt}}_{\sigma_2}$}

\end{picture}

\noindent
where $\sigma_1 \in \mathcal{NC}^2_{2k_1}$, $\sigma_2 \in \mathcal{NC}^2_{2k_2}$ and $k_1 + k_2 = pk-1$. The inverse of $\varphi$ is well defined, which means that $\varphi$ is a bijection on $\mathcal{NC}^2_{2pk}$. By $\varphi^n$ we denote the $n$-th composition of $\varphi$.

\begin{Lemma}
For any $i\in \{1,2, \ldots , p\}$, it holds that
$$
d_i (\mathcal{N}_i(x,{\bf d}) -1) = d_0(\mathcal{N}_0(x,{\bf d})-1).
$$ 
\end{Lemma}
{\it Proof.} First, we will show that 
\begin{equation}
N_i(k, j_0+1, j_1, \ldots, j_p) = N_0(k, j_0,\ldots, j_{i-1}, j_i+1,j_{i+1}, \ldots, j_p).
\end{equation}
Let $\pi \in \mathcal{NC}^2_{2pk}(W_i^k)$. Then the partition $\pi$ given by 
the diagram (the legs are labelled only by the set of letters $\{1,\ldots , p, p^{*}, \ldots , 1\}$
and not by all consecutive numbers from the set $[2kp]$)

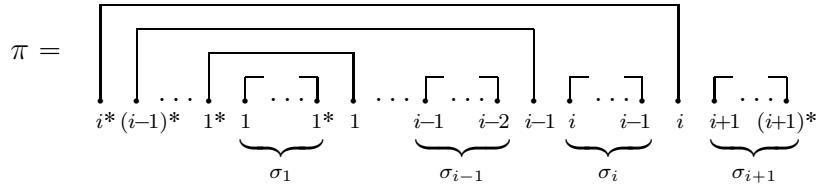
\unitlength=0.8mm
\special{em:linewidth 0.4pt}
\linethickness{0.4pt}
\begin{picture}(120.00,35.00)(-45.00,-3.00)

\put(-15,17){$\pi = $}

\put(00.00,10.00){\line(0,1){16}}
\put(00.00,26.00){\line(1,0){96}}
\put(96.00,10.00){\line(0,1){16}}
\put(00.00,10.00){\circle*{1.00}}
\put(96.00,10.00){\circle*{1.00}}
\put(-00.80,5.50){\tiny $i^*$}
\put(95.70,5.50){\tiny $i$}

\put(06.00,10.00){\line(0,1){12}}
\put(06.00,22.00){\line(1,0){66}}
\put(72.00,10.00){\line(0,1){12}}
\put(06.00,10.00){\circle*{1.00}}
\put(72.00,10.00){\circle*{1.00}}
\put(03.20,5.50){\tiny $(i\!\!-\!\!1)^*$}
\put(70.50,5.50){\tiny $i\!\!-\!\!1$}

\put(09.50,10.00){$\ldots$}

\put(18.00,10.00){\line(0,1){8}}
\put(18.00,18.00){\line(1,0){24}}
\put(42.00,10.00){\line(0,1){8}}
\put(18.00,10.00){\circle*{1.00}}
\put(42.00,10.00){\circle*{1.00}}
\put(17.00,5.50){\tiny $1^*$}
\put(41.40,5.50){\tiny $1$}

\put(24.00,10.00){\line(0,1){4}}
\put(24.00,14.00){\line(1,0){3}}
\put(36.00,10.00){\line(0,1){4}}
\put(36.00,14.00){\line(-1,0){3}}
\put(28.00,10.00){$\ldots$}
\put(24.00,10.00){\circle*{1.00}}
\put(36.00,10.00){\circle*{1.00}}
\put(23.20,5.50){\tiny $1$}
\put(35.00,5.50){\tiny $1^*$}
\put(22.80,04.00){$\underbrace{\hspace{33pt}}_{\sigma_1}$}

\put(45.50,10.00){$\ldots$}

\put(54.00,10.00){\line(0,1){4}}
\put(54.00,14.00){\line(1,0){3}}
\put(66.00,10.00){\line(0,1){4}}
\put(66.00,14.00){\line(-1,0){3}}
\put(58.00,10.00){$\ldots$}
\put(54.00,10.00){\circle*{1.00}}
\put(66.00,10.00){\circle*{1.00}}
\put(52.00,5.50){\tiny $i\!\!-\!\!1$}
\put(62.80,5.50){\tiny $i\!\!-\!\!2$}
\put(52.00,04.00){$\underbrace{\hspace{37pt}}_{\sigma_{i-1}}$}

\put(78.00,10.00){\line(0,1){4}}
\put(78.00,14.00){\line(1,0){3}}
\put(90.00,10.00){\line(0,1){4}}
\put(90.00,14.00){\line(-1,0){3}}
\put(82.00,10.00){$\ldots$}
\put(78.00,10.00){\circle*{1.00}}
\put(90.00,10.00){\circle*{1.00}}
\put(77.80,5.50){\tiny $i$}
\put(86.30,5.50){\tiny $i\!\!-\!\!1$}
\put(77.00,04.00){$\underbrace{\hspace{34pt}}_{\sigma_i}$}

\put(102.00,10.00){\line(0,1){4}}
\put(102.00,14.00){\line(1,0){3}}
\put(114.00,10.00){\line(0,1){4}}
\put(114.00,14.00){\line(-1,0){3}}
\put(106.00,10.00){$\ldots$}
\put(102.00,10.00){\circle*{1.00}}
\put(114.00,10.00){\circle*{1.00}}
\put(101.20,5.50){\tiny $i\!\!+\!\!1$}
\put(109.00,5.50){\tiny $(i\!\!+\!\!1)^*$}
\put(101.00,04.00){$\underbrace{\hspace{35pt}}_{\sigma_{i+1}}$}

\end{picture}

\noindent
is mapped by $\varphi^{\,i}$ onto 

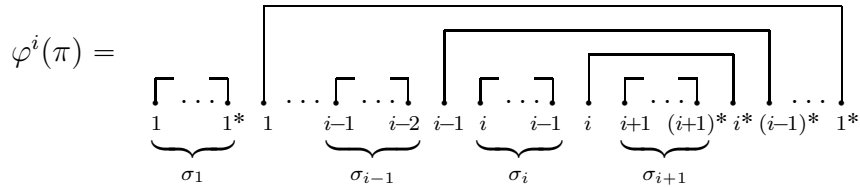
\begin{picture}(120.00,43.00)(-20.00,-8.00)

\put(0.00,17.00){$\varphi^i(\pi) = $}

\put(42.00,10.00){\line(0,1){16}}
\put(42.00,26.00){\line(1,0){96}}
\put(138.00,10.00){\line(0,1){16}}
\put(42.00,10.00){\circle*{1.00}}
\put(138.00,10.00){\circle*{1.00}}
\put(41.80,5.50){\tiny $1$}
\put(137.00,5.50){\tiny $1^*$}

\put(126.00,10.00){\line(0,1){12}}
\put(126.00,22.00){\line(-1,0){54}}
\put(72.00,10.00){\line(0,1){12}}
\put(126.00,10.00){\circle*{1.00}}
\put(72.00,10.00){\circle*{1.00}}
\put(124.20,5.50){\tiny $(i\!\!-\!\!1)^*$}
\put(70.20,5.50){\tiny $i\!\!-\!\!1$}

\put(129.50,10.00){$\ldots$}

\put(96.00,10.00){\line(0,1){8}}
\put(96.00,18.00){\line(1,0){24}}
\put(120.00,10.00){\line(0,1){8}}
\put(96.00,10.00){\circle*{1.00}}
\put(120.00,10.00){\circle*{1.00}}
\put(95.20,5.50){\tiny $i$}
\put(120.00,5.50){\tiny $i^*$}

\put(24.00,10.00){\line(0,1){4}}
\put(24.00,14.00){\line(1,0){3}}
\put(36.00,10.00){\line(0,1){4}}
\put(36.00,14.00){\line(-1,0){3}}
\put(28.00,10.00){$\ldots$}
\put(24.00,10.00){\circle*{1.00}}
\put(36.00,10.00){\circle*{1.00}}
\put(23.20,5.50){\tiny $1$}
\put(35.00,5.50){\tiny $1^*$}
\put(23.00,04.00){$\underbrace{\hspace{33pt}}_{\sigma_1}$}

\put(45.50,10.00){$\ldots$}

\put(54.00,10.00){\line(0,1){4}}
\put(54.00,14.00){\line(1,0){3}}
\put(66.00,10.00){\line(0,1){4}}
\put(66.00,14.00){\line(-1,0){3}}
\put(58.00,10.00){$\ldots$}
\put(54.00,10.00){\circle*{1.00}}
\put(66.00,10.00){\circle*{1.00}}
\put(52.00,5.50){\tiny $i\!\!-\!\!1$}
\put(62.80,5.50){\tiny $i\!\!-\!\!2$}
\put(52.00,04.00){$\underbrace{\hspace{37pt}}_{\sigma_{i-1}}$}

\put(78.00,10.00){\line(0,1){4}}
\put(78.00,14.00){\line(1,0){3}}
\put(90.00,10.00){\line(0,1){4}}
\put(90.00,14.00){\line(-1,0){3}}
\put(82.00,10.00){$\ldots$}
\put(78.00,10.00){\circle*{1.00}}
\put(90.00,10.00){\circle*{1.00}}
\put(77.80,5.50){\tiny $i$}
\put(86.30,5.50){\tiny $i\!\!-\!\!1$}
\put(77.00,04.00){$\underbrace{\hspace{34pt}}_{\sigma_i}$}

\put(102.00,10.00){\line(0,1){4}}
\put(102.00,14.00){\line(1,0){3}}
\put(114.00,10.00){\line(0,1){4}}
\put(114.00,14.00){\line(-1,0){3}}
\put(106.00,10.00){$\ldots$}
\put(102.00,10.00){\circle*{1.00}}
\put(114.00,10.00){\circle*{1.00}}
\put(101.20,5.50){\tiny $i\!\!+\!\!1$}
\put(109.00,5.50){\tiny $(i\!\!+\!\!1)^*$}
\put(101.00,04.00){$\underbrace{\hspace{35pt}}_{\sigma_{i+1}}$}

\end{picture}

\noindent
where $\sigma_r \in \mathcal{NC}^2_{2k_r}$ and $k_1 + \ldots + k_{i+1} = pk-i$. 
We understand that certain partitions among $\sigma_{1}, \ldots , \sigma_{i+1}$ in these diagrams may be empty. 
Obviously, $\varphi^{\,i}(\pi)$ is adapted to $W_0^k$, since $\pi$ is adapted to $W_i^k$ and, conversely, if $\pi$ wasn't adapted to 
$W_i^k$, then $\varphi^{\,i}(\pi)$ wouldn't be adapted to $W_0^k$. This means that $\varphi^{\,i}$ restricted to the set $\mathcal{NC}^2_{2pk}(W_i^k)$ 
is a bijection between $\mathcal{NC}^2_{2pk}(W_i^k)$ and $\mathcal{NC}^2_{2pk}(W_0^k)$. Moreover, as compared with $\pi$, 
$\varphi^i(\pi)$ has one more block labelled by the ordered pair $(r,r^*)$ and one less block labelled by $(r^*,r)$ for each $r\in\{1,\ldots,i\}$. 
The numbers of other blocks are the same in $\pi$ and $\varphi^{i}(\pi)$.
This means that 
$$
j_{r}(\varphi^{i}(\pi))=
\left\{
\begin{array}{ll}
j_0(\pi)-1& {\rm if}\;r=0\\
j_i(\pi)+1& {\rm if}\;r=i\\
j_i(\pi)  & {\rm otherwise}
\end{array}
\right.
$$ 
which proves (4). Now, using the fact that $j_0(\pi)\geqslant1$ for $\pi \in \mathcal{NC}^2_{2pk}(W_i^k)$ and $i,k\geqslant1$, we get
\begin{eqnarray*}
d_i(\mathcal{N}_i(x,{\bf d}) \!-\! 1) 
&=& 
d_i \sum_{k=1}^{\infty} \sum_{j_{0}=1}^{\infty}
\sum_{j_1, \ldots , j_p=0}^{\infty}
N_i(k,j_0,\ldots,j_p) d_0^{j_0} \ldots d_p^{j_p}x^{k} \\
&=& 
d_0 \sum_{k=1}^{\infty} 
\sum_{j_0, \ldots , j_p=1}^{\infty}N_i(k,j_0+1,\ldots,j_p)\\
&&
\times d_0^{j_0} \ldots d_{i-1}^{j_{i-1}} d_{i}^{j_i+1} d_{i+1}^{j_{i+1}} \ldots d_p^{j_p}x^{k}
\end{eqnarray*}
Here we use (4) and the observation that $j_i(\pi)\geqslant1$ for $\pi \in \mathcal{NC}^2_{2pk}(W_0^k)$ and $i,k\geqslant1$.
This becomes 
\begin{eqnarray*}
d_i (\mathcal{N}_i(x,{\bf d}) \!-\! 1) 
&=&
d_0 \sum_{k=1}^{\infty} \sum_{j_0,\ldots,j_p=0}^{\infty}N_0(k,j_0,\ldots,j_{i},j_i+1,\ldots,j_p) \\
&& \times 
d_0^{j_0} \ldots  d_i^{j_i+1}  \ldots d_p^{j_p}x^{k} \\
&=&
d_0 \sum_{k=1}^{\infty}  \sum_{j_i=1}^{\infty} \ \sum_{j_0,\ldots,j_{i-1},j_{i+1}\ldots,j_p=0}^{\infty} \!\! 
N_0(k,j_0,\ldots,j_i,\ldots,j_p) \\
&& \times
d_0^{j_0} \ldots  d_i^{j_i}  \ldots d_p^{j_p}x^{k} \\
&=&
d_0 \sum_{k=1}^{\infty} \sum_{j_0,\ldots,j_p=0}^{\infty} N_0(k,{\bf j}) {\bf d^j}x^{k} \\
&=& d_0 (\mathcal{N}_0(x,{\bf d}) -1),
\end{eqnarray*} 
which completes the proof. \hfill $\blacksquare$\\

\begin{Lemma}
For any $p\in {\mathbb N}$, the generating functions ${\mathcal N}_{0}, \ldots , {\mathcal N}_{p}$ satisfy the equation
$$
\mathcal{N}_0(x,{\bf d})-1 = xd_1\ldots d_p \, \prod_{i=0}^{p} \mathcal{N}_i(x,{\bf d}).
$$
\end{Lemma}
{\it Proof.}
In order to prove this lemma, let us establish a recurrence relation between numbers $N_i(k,{\bf j})$ for 
$i\in \{0, \ldots , p\}$ and certain values of $k$ and ${\bf j}=(j_0,\ldots , j_p)$. 
For that purpose, observe that each partition $\pi \in \mathcal{NC}^2_{2pk}(W_0^k)$ can be expressed in terms of the diagram

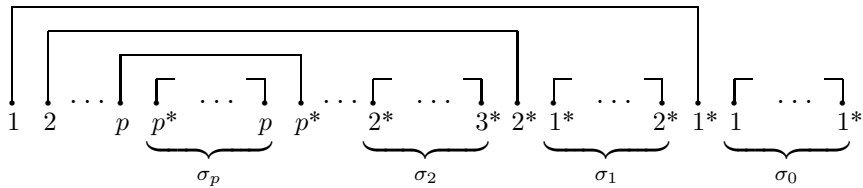
\unitlength=0.8mm
\special{em:linewidth 0.4pt}
\linethickness{0.4pt}
\begin{picture}(120.00,42.00)(-20.00,-7.00)

\put(00.00,10.00){\line(0,1){16}}
\put(00.00,26.00){\line(1,0){114}}
\put(114.00,10.00){\line(0,1){16}}
\put(00.00,10.00){\circle*{1.00}}
\put(114.00,10.00){\circle*{1.00}}
\put(-00.80,5.50){\footnotesize $1$}
\put(113.00,5.50){\footnotesize $1^*$}

\put(06.00,10.00){\line(0,1){12}}
\put(06.00,22.00){\line(1,0){78}}
\put(84.00,10.00){\line(0,1){12}}
\put(06.00,10.00){\circle*{1.00}}
\put(84.00,10.00){\circle*{1.00}}
\put(05.20,5.50){\footnotesize $2$}
\put(83.00,5.50){\footnotesize $2^*$}

\put(09.50,10.00){$\ldots$}

\put(18.00,10.00){\line(0,1){8}}
\put(18.00,18.00){\line(1,0){30}}
\put(48.00,10.00){\line(0,1){8}}
\put(18.00,10.00){\circle*{1.00}}
\put(48.00,10.00){\circle*{1.00}}
\put(17.20,5.50){\footnotesize $p$}
\put(47.00,5.50){\footnotesize $p^*$}

\put(24.00,10.00){\line(0,1){4}}
\put(24.00,14.00){\line(1,0){3}}
\put(42.00,10.00){\line(0,1){4}}
\put(42.00,14.00){\line(-1,0){3}}
\put(31.00,10.00){$\ldots$}
\put(24.00,10.00){\circle*{1.00}}
\put(42.00,10.00){\circle*{1.00}}
\put(23.20,5.50){\footnotesize $p^*$}
\put(41.00,5.50){\footnotesize $p$}
\put(22.00,04.00){$\underbrace{\hspace{49pt}}_{\sigma_p}$}

\put(51.50,10.00){$\ldots$}

\put(60.00,10.00){\line(0,1){4}}
\put(60.00,14.00){\line(1,0){3}}
\put(78.00,10.00){\line(0,1){4}}
\put(78.00,14.00){\line(-1,0){3}}
\put(67.00,10.00){$\ldots$}
\put(60.00,10.00){\circle*{1.00}}
\put(78.00,10.00){\circle*{1.00}}
\put(59.20,5.50){\footnotesize $2^*$}
\put(77.00,5.50){\footnotesize $3^*$}
\put(58.00,04.00){$\underbrace{\hspace{49pt}}_{\sigma_2}$}

\put(90.00,10.00){\line(0,1){4}}
\put(90.00,14.00){\line(1,0){3}}
\put(108.00,10.00){\line(0,1){4}}
\put(108.00,14.00){\line(-1,0){3}}
\put(97.00,10.00){$\ldots$}
\put(90.00,10.00){\circle*{1.00}}
\put(108.00,10.00){\circle*{1.00}}
\put(89.20,5.50){\footnotesize $1^*$}
\put(106.50,5.50){\footnotesize $2^*$}
\put(88.00,04.00){$\underbrace{\hspace{49pt}}_{\sigma_1}$}

\put(120.00,10.00){\line(0,1){4}}
\put(120.00,14.00){\line(1,0){3}}
\put(138.00,10.00){\line(0,1){4}}
\put(138.00,14.00){\line(-1,0){3}}
\put(127.00,10.00){$\ldots$}
\put(120.00,10.00){\circle*{1.00}}
\put(138.00,10.00){\circle*{1.00}}
\put(119.20,5.50){\footnotesize $1$}
\put(137.00,5.50){\footnotesize $1^*$}
\put(118.00,04.00){$\underbrace{\hspace{49pt}}_{\sigma_0}$}

\end{picture}

\vspace{0pt}
\noindent
where $\sigma_i \in \mathcal{NC}^2_{2pk_i}(W_i^{k_i})$ and $k_0 + k_1 + \ldots + k_p = k-1$. We understand that each partition
$\sigma_{i}$ in this diagram may be empty. Let $j_{i,r}=j_{r}(\sigma_{i})$ be the number of blocks labelled by $(r,r^{*})$ or 
$((r+1)^{*},r+1)$ in $\sigma_{i}$. It will be convenient to denote by
$$
S_r(\pi) = j_{0,r}+j_{1,r}+\ldots +j_{p,r}
$$ 
where $r\in \{0,1, \ldots , p\}$, the numbers of blocks labelled by $(r,r^{*})$ or $((r+1)^{*}, r+1)$, respectively, 
in the union $\sigma_{0}\cup \ldots \cup \sigma_{p}$. Note that 
$$
S_{r}(\pi)=\left\{
\begin{array}{ll}
j_{r}(\pi)-1& {\rm if}\;r\geqslant1\\
j_{0}(\pi) & {\rm if}\; r=0
\end{array}
\right.,
$$
which leads to the recurrence formula
\begin{equation}
N_0(k,{\bf j}) = \sum_{k_0 + \ldots + k_p = k-1} 
\sum_{\stackrel{{\bf j}_{0}, \ldots , {\bf j}_{p}}
{\scriptscriptstyle S_0 = j_0,S_1 = j_1 - 1,  \ldots ,S_p = j_p - 1} }
\prod_{i=0}^{p}N_i(k_i,{\bf j}_{\,i})
\end{equation}
for $k\geqslant 1$, where ${\bf j}_{\,i}=(j_{i,0},\ldots,j_{i,p})$ and 
$S_r = j_{0,r}+j_{1,r}+\ldots +j_{p,r}$. Using this formula, we get
\begin{eqnarray*}
\mathcal{N}_0(x,{\bf d}) \!-\!1 \!\!\!&=&\!\! \sum_{k=0}^{\infty}\sum_{{\bf j}}
N_0(k+1,{\bf j}) {\bf d^j}x^{k+1} \\
&=&
\sum_{k=0}^{\infty}
\sum_{k_0+\ldots + k_p=k}
\sum_{\stackrel{{\bf j},{\bf j}_{0}, \ldots , {\bf j}_{p}}
{\scriptscriptstyle S_0 = j_0,S_1 = j_1 - 1,  \ldots , S_p = j_p - 1} }
\prod_{i=0}^{p}N_i(k_i,{\bf j}_{\,i}){\bf d}^{{\bf j}}x^{k+1} \\
&=&
x d_1 d_2 \ldots d_p
\sum_{k=0}^{\infty}
\sum_{k_0 + \ldots + k_p = k} 
\sum_{\stackrel{{\bf j},{\bf j}_{0}, \ldots , {\bf j}_{p}}
{\scriptscriptstyle S_0 = j_0,S_1 = j_1,  \ldots , S_p = j_p} }
\prod_{i=0}^{p}N_i(k_i,{\bf j}_{\,i}){\bf d}^{\bf j}x^{k} \\
&=&
x d_1d_2 \ldots d_p  \prod_{i=0}^{p}\mathcal{N}_i(x,{\bf d}) 
\end{eqnarray*}
which proves our assertion.
\hfill $\blacksquare$\\

\section{Main results}

The main result of this paper consists in demonstrating that the multivariate Fuss-Narayana polynomials $F_k$
defined in Section 2 are special cases of the polynomials $P_k$ obtained in [8] 
as limit moments of products of independent random matrices. This fact is proved in the theorem given below. We
also show that the polynomials $F_k$ are moments of free multiplicative convolutions of 
Marchenko-Pastur distributions with arbitrary shape parameters.

\begin{Theorem}
The coefficients of the polynomial $P_k$ are generalized Fuss-Narayana numbers, namely
$$
N_0(k,{\bf j}) = \frac{1}{k} 
{k \choose j_0+1}{k \choose j_1}\ldots {k \choose j_{p}}
$$
where ${\bf j}=(j_0, j_1, \ldots , j_p)$, whenever $k\geqslant 1$ and $j_0 + \ldots + j_p = pk$, 
with $N_0(0,0,\ldots,0) = 1$. In other cases, $N_0(k,{\bf j}) = 0$.
\end{Theorem}

{\it Proof.} By Lemmas 3.1 and 3.2, we have
$$
d_0 (\mathcal{N}_0(x,{\bf d}) - 1) = x \prod_{i=0}^{p}\left( d_i (\mathcal{N}_0(x,{\bf d}) - 1) + d_i \right),
$$
which means that the function of the form
$$
g(x)=d_0 (\mathcal{N}_0(x,{\bf d}) - 1)
$$ 
is a solution of the equation (2). Hence, by Proposition 2.2, 
$$
d_0 (\mathcal{N}_0(x,{\bf d}) - 1) = \mathcal{N}(x,{\bf d}).
$$
Thus,
\begin{eqnarray*}
\sum_{k=1}^{\infty} \sum_{j_0,\ldots,j_{p}=0}^{\infty} N_0(k,{\bf j}) {\bf d^j}x^{k} & =& 
\sum_{k=1}^{\infty} \ \sum_{j_0 + \ldots + j_{p} = pk+1} 
N(k,j_0,\ldots , j_{p})d_0^{j_0 - 1} d_1^{j_1} \ldots d_{p}^{j_{p}}x^{k}\\
&=& \sum_{k=1}^{\infty} \ \sum_{j_0 + \ldots + j_{p} = pk} N(k,j_0+1, \ldots, j_{p}) {\bf d^j}x^{k}
\end{eqnarray*}
since the generalized Fuss-Narayana numbers are defined for $j_0,\ldots , j_{p}\in {\mathbb N}$ and thus one can divide 
${\mathcal N}(x,{\bf y})$ by $y_{0}$ and then replace the summation index $j_{0}$ by $j_{0}+1$. Therefore, in the last
summation, $j_0\in {\mathbb N}\cup \{0\}$. Of course, also $j_i\in {\mathbb N}\cup \{0\}$ for $i>0$, with the restriction that
$j_0+\ldots +j_{p}=kp$ and $j_0, \ldots , j_{p}\leq k$, which means that only one index among these may be equal to zero. 
Comparing the coefficients we obtain our assertion for $k\geqslant1$. The case $k=0$ is obvious. \hfill $\blacksquare$\\

\begin{Corollary}
For any $k,p\in {\mathbb N}$, it holds that 
$$
F_{k}(t_1,t_2, \ldots , t_p)=P_{k}(1,t_1, \ldots , t_p)
$$
\end{Corollary}
{\it Proof.}
By Theorem 4.1, $d_{0}P_{k}(d_0, d_1, \ldots , d_p)=R_{k}(d_0, d_1, \ldots , d_p)$ and thus, using 
Remark 2.1, we get the desired formula.
\hfill $\blacksquare$\\

\begin{Example}
{\rm In Example 3.1, the polynomials $P_2$ and $P_3$ were calculated for $p=2$ by means of Definition 3.2.
Now, we can calculate them using Theorem 4.1. We have
$$
P_k(d_0,d_1,d_2)=\sum_{j_0+j_1+j_2=pk}N_{0}(k, j_0,j_1,j_2)d_{0}^{j_{0}}d_{1}^{j_2}d_{2}^{j_3}.
$$
For $k=2$, we obtain $N_0(2,0,2,2)=N_0(2,1,1,2)=N_0(2,1,2,1)=1$ and thus
$$
F_{2}(t_1,t_2)=t_1^2t_2^2+t_1t_2^2+t_1^2t_2
$$
is the corresponding multivariate Fuss-Narayana polynomial.
In turn, for $k=3$, we obtain $N_0(3,0,3,3)=N_0(3,2,1,3)=N_0(3,2,3,1)=1$ and 
$N_0(3,2,2,2)=N_0(3,1,2,3)=N_0(3,1,3,2)=3$ and thus
$$
F_{3}(t_1,t_2)=t_1^3 t_2^3 + t_1 t_2^3 + t_1^3 t_2 + 3 t_1^2 t_2^2 + 3 t_1^2 t_2^3 +3 t_1^3 t_2^2
$$
is the corresponding mutlivariate Fuss-Narayana polynomial.
}
\end{Example}

The next corollary is an easy application of the above theorem and of the Vandermonde's identity. The original proof
of this fact is due to Kemp and Speicher [6].

\begin{Corollary}
The cardinality of the set $\mathcal{NC}^2_{2pk}(W_0^k)$ is the Fuss-Catalan number, i.e.
$$
|\mathcal{NC}^2_{2pk}(W_0^k)| = \frac{1}{pk+1} \left( \genfrac{}{}{0pt}{}{(p+1)k}{k} \right).
$$
\end{Corollary}
{\it Proof.} It is enough to put $d_0 = d_1 = \ldots = d_p = 1$ in (3). \hfill $\blacksquare$\\

Finally, let us show a direct application of multivariate Fuss-Catalan polynomials to free probability. Denote by
$$
\varrho_{t}={\rm max}\{1-t,0\}\delta_{0}+\frac{\sqrt{(x-a)(b-x)}}{2\pi x}1\!\!1_{[a,b]}(x)dx
$$
where $a=(1-\sqrt{t})^{2}$ and $b=(1+\sqrt{t})^{2}$, the Marchenko-Pastur law with the 
shape parameter equal to $t>0$. This distribution plays the role of the free analog of the Poisson law and is often called 
the free Poisson law [11].

We will prove that the moments of free multiplicative convolutions 
$$
\varrho_{t_1}\boxtimes \varrho_{t_2}\boxtimes \ldots \boxtimes \varrho_{t_p}
$$
of Marchenko-Pastur laws with different shape parameters $t_1,t_2,\ldots , t_p$ are multivariate Fuss-Narayana polynomials. 
An explicit formula for the much simpler case when all shape parameters are equal (and thus we deal with convolution powers of $\varrho_t$) 
has recently been found by Hinz and M{\l}otkowski [5]. Another special case of our formula is that of Fuss-Narayana polynomials, 
obtained for $t_1=\ldots =t_{p-1}=1$ and $t_p=t$, which gives the moments of free Bessels laws $\pi_{p,t}$ of 
Banica {\it et al} [3, Theorem 5.2].

\begin{Proposition}
For any positive $t_1,t_2, \ldots, t_p$ and $k\in {\mathbb N}$, it holds that
$$
m_{k}(\varrho_{t_1}\boxtimes \varrho_{t_2}\boxtimes \ldots \boxtimes \varrho_{t_{p}})=F_{k}(t_1,t_2, \ldots , t_p).
$$
\end{Proposition}
{\it Proof.}
Let $S_{\mu}$ and $\psi_{\mu}$ denote the S-transform and the moment generating function (without the constant term) of a probability measure on the real line $\mu$, respectively. It is well-known that 
$$
S_{\mu}(z)=\frac{1+z}{z}\psi_{\mu}^{-1}(z),
$$
where $\psi_{\mu}^{-1}$ is the composition inverse of $\psi_{\mu}$.
Moreover, the S-transform of the Marchenko-Pastur law is 
$$
S_{\varrho_{t}}(z)=\frac{1}{z+t}
$$
for any $t>0$. Using the mutliplicativity of the S-transform with respect to the free multiplicative convolution, we obtain
$$
\psi_{\varrho}^{-1}(z)=\frac{1}{(z+1)(z+t_1)\ldots (z+t_p)}
$$
where $\varrho=\varrho_{t_1}\boxtimes \varrho_{t_2}\boxtimes \ldots \boxtimes \varrho_{t_{p}}$. Obviously, the above formula is equivalent to
$$
\psi_{\varrho}(z)=
z(\psi_{\varrho}(z)+1)(\psi_{\varrho}(z)+t_1)\ldots (\psi_{\varrho}(z)+t_p).
$$
Now, using Proposition 2.2 with $d_{0}=1$, $d_1=t_1, \ldots , d_{p}=t_{p}$,
we obtain the desired formula.\hfill $\blacksquare$\\

\begin{Example}
{\rm $m_2(\varrho_{t_1}\boxtimes \varrho_{t_2}\boxtimes \varrho_{t_{3}})=
t_1^2t_2^2t_3^2+t_1t_2^2t_3^2+t_1^2t_2t_3^2+t_1^2t_2^2t_3$.}
\end{Example}

\end{document}